\pdfoutput=1

\documentclass{article}
\usepackage[utf8]{inputenc}
\usepackage[T2A]{fontenc}
\usepackage[english]{babel}
\usepackage[intlimits]{amsmath}
\usepackage{amssymb}
\usepackage[numbers]{natbib}
\usepackage[ruled, linesnumbered]{algorithm2e}
\usepackage{graphicx}

\newtheorem{definition}{Definition}
\newenvironment{thm}{\begin{theorem}}{\end{theorem}}
\newtheorem{theorem}{Theorem}

\newtheorem{corollary}{Corollary}

\DeclareMathOperator*{\diag}{\mathrm{diag}\,}
\DeclareMathOperator*{\vspan}{\mathrm{span}\,}
\DeclareMathOperator*{\dist}{\mathrm{dist}\,}

\begin{document}
\title{Computation- and Space-Efficient Implementation of SSA}
\author{Anton Korobeynikov \\
    Department of Statistical Modelling, \\
    Saint Petersburg State University \\
    198504, Universitetskiy pr. 28, Saint Petersburg, Russia\\
    \texttt{asl@math.spbu.ru}}

\maketitle

\begin{abstract}
    The computational complexity of different steps of the basic SSA is
    discussed. It is shown that the use of the general-purpose ``blackbox''
    routines (e.g. found in packages like \verb@LAPACK@) leads to huge waste of
    time resources since the special Hankel structure of the trajectory matrix
    is not taken into account.

    We outline several state-of-the-art algorithms (for example, Lanczos-based
    truncated SVD) which can be modified to exploit the structure of the
    trajectory matrix. The key components here are Hankel matrix-vector
    multiplication and the Hankelization operator. We show that both can be
    computed efficiently by the means of Fast Fourier Transform.

    The use of these methods yields the reduction of the worst-case
    computational complexity from $O(N^{3})$ to $O(k N \log{N})$, where $N$ is
    series length and $k$ is the number of eigentriples desired.
\end{abstract}

\section{Introduction}
\label{sec:introduction}
Despite the recent growth of SSA-related techniques it seems that little
attention was paid to their efficient implementations, e.g. the singular value
decomposition of Hankel trajectory matrix, which is a dominating time-consuming
step of the basic SSA.

In almost all SSA-related applications only a few leading eigentriples are
desired, thus the use of general-purpose ``blackbox'' SVD implementations causes
a huge waste of time and memory resources. This almost prevents the use of the
optimum window sizes even on moderate length (say, few thousand) time
series. The problem is much more severe for 2D-SSA \cite{Golyandina09a} or
subspace-based methods \cite{Badeau08, Djermoune09, Golyandina09b}, where large
window sizes are typical.

We show that proper use of state-of-art singular value decomposition algorithms
can significantly reduce the amount of operations needed to compute truncated
SVD suitable for SSA purposes. This is possible because the information about
the leading singular triples becomes available long before the costly full
decomposition.

These algorithms are usually based on Lanczos iterations or related
methods. This algorithm can be modified to exploit the Hankel structure of the
data matrix effectively. Thus the computational burden of SVD can be
considerably reduced.

In section~\ref{sec:complexity} we will outline the computational and storage
complexity of the steps of the basic SSA algorithm. Section~\ref{sec:svd}
contains a review of well-known results connected with the Lanczos singular
value decompostion algorithm. Generic implementation issues in finite-precision
arithmetics will be discussed as well. The derivation of the efficient Hankel
SSA trajectory matrix-vector multiplication method will be presented in
section~\ref{sec:matmul}. This method is a key component of the fast truncated
SVD for the Hankel SSA trajectory matrices. In the section~\ref{sec:diagavg} we
will exploit the special structure of the rank 1 Hankelization operator which
allows its efficient implementation. Finally, all the implementations of the
mentioned algorithms are compared in sections~\ref{sec:impl}
and~\ref{sec:example} on artificial and real data.

\section{Complexity of the Basic SSA Algorithm}
\label{sec:complexity}
Let $N > 2$. Denote by $F = \left(f_{1},\dots,f_{N}\right)$ a real-valued time
series of length $N$. Fix the \emph{window length} $L$, $1 < L < N$.  The Basic
SSA algorithm for the decomposition of a time series consists of four steps
\cite{Golyandina01}. We will consider each step separately and discuss its
computational and storage complexity.

\subsection{Embedding}
The embedding step maps original time series to a sequence of lagged
vectors. The embedding procedure forms $K = N - L + 1$ \emph{lagged vectors}
$X_{i} \in \mathbb{R}^{L}$:
\begin{equation*}
    X_{i} = \left(f_{i},\dots,f_{i+L-1}\right)^{T}, \quad 1 \leqslant i
    \leqslant K.
\end{equation*}
These vectors form the \emph{$L$-trajectory} (or \emph{trajectory}) matrix $X$
of the series $F$:
\begin{equation*}
    X =
    \begin{pmatrix}
        f_{1} & f_{2} & \cdots & f_{K-1} & f_{K}   \\
        f_{2} & f_{3} & \cdots & f_{K}  & f_{K+1} \\
        \vdots & \vdots & & \vdots  & \vdots \\
        f_{L} & f_{L+1} & \cdots & f_{N-1} & f_{N}
    \end{pmatrix}.
\end{equation*}
Note that this matrix has equal values along the antidiagonals, thus it is a
\emph{Hankel matrix}. The computational complexity of this step is negligible:
it consists of simple data movement which is usually considered a cheap
procedure. However, if the matrix is stored as-is, then this step obviously has
storage complexity of $O(LK)$ with the worst case being $O(N^{2})$ when $K \sim
L \sim N/2$.

\subsection{Singular Value Decomposition}
This step performs the singular value decomposition of the $L \times K$
trajectory matrix $X$.  Without loss of generality we will assume that $L
\leqslant K$. Then the SVD can be written as
\begin{equation*}
    X = U^{T} \Sigma V,
\end{equation*}
where $U = \left(u_{1},\dots,u_{L}\right)$ is a $L \times L$ orthogonal matrix,
$V = \left(v_{1},\dots,v_{K}\right)$ is a $K \times K$ orthogonal matrix, and
$\Sigma$ is a $L \times K$ diagonal matrix with nonnegative real diagonal
entries $\Sigma_{ii} = \sigma_{i}$, for $i = 1, \dots, L$. The vectors $u_{i}$
are called \emph{left singular vectors}, the $v_{i}$ are the \emph{right
    singular vectors}, and the $\sigma_{i}$ are the \emph{singular values}. We
will label the singular values in descending order: $\sigma_{1} \geqslant
\sigma_{2} \geqslant \cdots \geqslant \sigma_{L}$. Note that the singular value
decomposition of $X$ can be represented in the form
\begin{equation*}
    X = X_{1} + \cdots + X_{L},
\end{equation*}
where $X_{i} = \sigma_{i}u_{i}v_{i}^{T}$. The matrices $X_{i}$ have rank $1$,
therefore we will call them \emph{elementary matrices}. The triple
$\left(\sigma_{i}, u_{i}, v_{i}\right)$ will be called a \emph{singular triple}.

Usually the SVD is computed by means of Golub and Reinsh algorithm
\cite{Golub70} which requires $O(L^{2}K + LK^{2} + K^{3})$ multiplications
having the worst computational complexity of $O(N^{3})$ when $K \sim L \sim
N/2$. Note that all the data must be processed at once and SVDs even of moderate
length time series (say, when $N > 5000$) are essentially unfeasible.

\subsection{Grouping}
The grouping procedure partitions the set of indices $\left\{1,\dots,L\right\}$
into $m$ disjoint sets $I_{1}, \dots, I_{m}$. For the index set $I =
\left\{i_{1}, \dots, i_{p}\right\}$ one can define the \emph{resultant matrix}
$X_{I} = X_{i_{1}} + \cdots + X_{i_{p}}$. Such matrices are computed for $I =
I_{1}, \dots, I_{m}$ and we obtain the decomposition
\begin{equation*}
    X = X_{I_{1}} + \cdots + X_{I_{m}}.
\end{equation*}

For the sake of simplicity we will assume that $m = L$ and each $I_{j} =
\left\{j\right\}$ later on. The computation of each elementary matrix $X_{j}$
costs $O(LK)$ multiplications and $O(LK)$ storage yielding overall computation
complexity of $O(L^{2}K)$ and $O(L^{2}K)$ for storage. Again the worst case is
achieved at $L \sim K \sim N/2$ with the value $O(N^{3})$ for both computation
and storage.

\subsection{Diagonal Averaging (Hankelization)}
\label{ssec:hankelization}
This is the last step of the Basic SSA algorithm transforming each grouped
matrix $X_{I_{j}}$ into a new time series of length $N$. This is performed by
means of \emph{diagonal averaging procedure} (or \emph{Hankelization}). It
transforms an arbitrary $L \times K$ matrix $Y$ to the series
$g_{1},\dots,g_{N}$ by the formula:
\begin{equation*}
    g_{k} =
    \begin{cases}
        \frac{1}{k}\sum_{j=1}^{k}{y_{j,k-j+1}}, & \quad 1 \leqslant k \leqslant L, \\
        \frac{1}{L}\sum_{j=1}^{L}{y_{j,k-j+1}}, & \quad L < k < K, \\
        \frac{1}{N-k+1}\sum_{j=k-K+1}^{L}{y_{j,k-j+1}}, & \quad K \leqslant k
        \leqslant N.
    \end{cases}
\end{equation*}

Note that these two steps are usually combined. One just forms the elementary
matrices one by one and immediately applies the Hankelization operator. This
trick reduces the storage requirements considerably.

$O(N)$ multiplications and $O(LK)$ additions are required for performing the
diagonal averaging of one elementary matrix. Summing these values for $L$
elementary matrices we obtain the overall computation complexity for this step
being $O(L^{2}K)$ additions and $O(NL)$ multiplications. The worst case
complexity will be then $O(N^{3})$ additions and $O(N^{2})$ multiplications

Summing the complexity values alltogether for all four steps we obtain the worst
case computational complexity to be $O(N^{3})$. The worst case occurs exactly
when window length $L$ equals $N/2$ (the optimal window length for
asymptotic separability \cite{Golyandina01}).

\section{Truncated Singular Value Decomposition}
\label{sec:svd}
The typical approach to the problem of computating the singular triples
$\left(\sigma_{i},u_{i},v_{i}\right)$ of $A$ is to use the Schur decomposition
of some matrices related to $A$, namely
\begin{enumerate}
    \item the \emph{cross product} matrices $A A^{T}$, $A^{T}A$, or
    \item the \emph{cyclic} matrix $C = \left(\begin{smallmatrix}0 & A \\ A^{T} &
            0 \end{smallmatrix}\right)$.
\end{enumerate}
One can prove the following theorem (see, e.g. \cite{Golub96, Larsen98}):
\begin{thm}
    \label{thm:conn}
    Let $A$ be an $L \times K$ matrix and assume without loss of generality that
    $K \geqslant L$. Let the singular value decomposition of $A$ be
    \begin{equation*}
        U^{T}AV = \diag\left(\sigma_{1},\dots,\sigma_{L}\right).
    \end{equation*}
    Then
    \begin{eqnarray*}
        V^{T}\left(A^{T}A\right)V & = &
        \diag(\sigma_{1}^{2},\dots,\sigma_{L}^{2},\underbrace{0,\dots,0}_{K-L}), \\
        U^{T}\left(AA^{T}\right)U & = &
        \diag(\sigma_{1}^{2},\dots,\sigma_{L}^{2}). \\
    \end{eqnarray*}
    Moreover, if $V$ is partitioned as
    \begin{equation*}
        V = \left(V_{1}, V_{2}\right),
        \quad V_{1} \in \mathbb{R}^{K \times L}, V_{2} \in \mathbb{R}^{K \times (K-L)},
    \end{equation*}
    then the orthonormal columns of the $(L+K) \times (L+K)$ matrix
    \begin{equation*}
        Y =
        \frac{1}{\sqrt{2}}
        \begin{pmatrix}
            U    &  - U   & 0 \\
            V_{1} & V_{1}  & \sqrt{2} V_{2}
        \end{pmatrix}
    \end{equation*}
    form an eigenvector basis for the cyclic matrix $C$ and
    \begin{equation*}
        Y^{T}CY =
        \diag
        (\sigma_{1},\dots,\sigma_{L},
        -\sigma_{1},\dots,-\sigma_{L},
        \underbrace{0,\dots,0}_{K-L}
        ).
    \end{equation*}
\end{thm}

The theorem tells us that we can obtain the singular values and vectors of $A$
by computing the eigenvalues and corresponding eigenvectors of one of the
symmetric matrices. This fact forms the basis of any $SVD$ algorithm.

In order to find all eigenvalues and eigenvectors of the real symmetric matrix
one usually turns it into a similar tridiagonal matrix by the means of
orthogonal transformations (Householder rotations or alternatively via Lanczos
recurrenses). Then, for example, another series of orthogonal transformations
can be applied to the latter tridiagonal matrix which converges to a diagonal
matrix \cite{Golub96}.

Such approach, which in its original form is due to Golub and Reinsh
\cite{Golub70} and is used in e.g. \verb@LAPACK@ \cite{Lapack99}, computes the
SVD by implicitly applying the QR algorithm to the symmetric eigenvalue problem
for $A^{T}A$.

For large and sparse matrices the Golub-Reinsh algorithm is impractical. The
algorithm itself starts out by applying a series of transformations directly to
matrix $A$ to reduce it to the special bidiagonal form. Therefore it requires
the matrix to be stored explicitlly, which may be impossible simply due to its
size. Moreover, it may be difficult to take advantage of any structure
(e.g. Hankel in SSA case) or sparsity in $A$ since it is quickly destroyed by
the transformations applied to the matrix. The same argument applies to SVD
computation via the Jacobi algorithm.

Bidiagonalization was proposed by Golub and Kahan \cite{Golub65} as a way of
tridiagonalizing the cross product matrix without forming it
explicitly. The method yields the decomposition
\begin{equation*}
    A = P^{T}BQ,
\end{equation*}
where $P$ and $Q$ are orthogonal matrices and $B$ is an $L \times K$ lower
bidiagonal matrix. Here $BB^{T}$ is similar to $AA^{T}$. Additionally, there are
well-known SVD algorithms for real bidiagonal matrices, for example, the QR
method \cite{Golub96}, divide-and-conquer \cite{Gu95}, and twisted
\cite{Dhillon04} methods.

The structure of our exposition follows \cite{Larsen98}.

\subsection{Lanczos Bidiagonalization}
\label{sec:lbidiag}
For a rectangular $L \times K$ matrix $A$ the Lanczos bidiagonalization computes
a sequence of \emph{left Lanczos vectors} $u_{j} \in \mathbb{R}^{L}$ and
\emph{right Lanczos vectors} $v_{j} \in \mathbb{R}^{K}$ and scalars $\alpha_{j}$
and $\beta_{j}$, $j=1,\dots,k$ as follows:

\begin{algorithm}[H]
    \caption{Golub-Kahan-Lanczos Bidiagonalization}
    \label{alg:gklbdiag}
    \dontprintsemicolon
    Choose a starting vector $p_{0} \in \mathbb{R}^{L}$ \;
    $\beta_{1} \longleftarrow \left\|p_{0}\right\|_{2}$, $u_{1} \longleftarrow
    p_{0} / \beta_{1}$, $v_{0} \longleftarrow 0$ \;
    \For{$j \leftarrow 1$ \KwTo $k$}{\nllabel{laniter}
        $r_{j} \longleftarrow A^{T}u_{j} - \beta_{j}v_{j-1}$ \;
        $\alpha_{j} \longleftarrow \left\|r_{j}\right\|_{2}$, \, $v_{j} \longleftarrow r_{j} / \alpha_{j}$ \;
        $p_{j} \longleftarrow Av_{j} - \alpha_{j}u_{j}$ \;
        $\beta_{j+1} \longleftarrow \left\|p_{j}\right\|_{2}$, \, $u_{j+1} \longleftarrow p_{j} / \beta_{j+1}$ \;
    }
\end{algorithm}
Later on we will refer to one iteration of the \textbf{for} loop at line
\ref{laniter} as a \emph{Lanczos step} (or \emph{iteration}). After $k$ steps we
will have the lower bidiagonal matrix
\begin{equation*}
    B_{k} =
    \begin{pmatrix}
        \alpha_{1} &            & & \\
        \beta_{2}  & \alpha_{2} &  & \\
                   & \beta_{3}  &  \ddots & \\
                   &            & \ddots  & \alpha_{k} \\
                   &            &         &  \beta_{k+1}
    \end{pmatrix}.
\end{equation*}
If all the computations performed are exact, then the Lanczos vectors will be
orthornomal:
\begin{equation*}
    U_{k+1} = \left(u_{1},\dots,u_{k+1}\right) \in \mathbb{R}^{L \times (k+1)},
    \quad
    U_{k+1}^{T}U_{k+1} = I_{k+1},
\end{equation*}
and
\begin{equation*}
    V_{k+1} = \left(v_{1},\dots,v_{k+1}\right) \in \mathbb{R}^{K \times (k+1)},
    \quad
    V_{k+1}^{T}V_{k+1} = I_{k+1},
\end{equation*}
where $I_{k}$ is $k \times k$ identity matrix. By construction, the columns of
$U_{k+1}$ and $V_{k+1}$ satisfy the recurrences
\begin{eqnarray*}
    \alpha_{j} v_{j} & = & A^{T}u_{j} - \beta_{j}v_{j-1}, \\
    \beta_{j+1} u_{j+1} & = & Av_{j} -\alpha_{j} u_{j},
\end{eqnarray*}
and can be written in the compact form as follows:
\begin{eqnarray}
    \label{eqn:lbdiagm1}
    AV_{k} & = & U_{k+1}B_{k}, \\
    \label{eqn:lbdiagm2}
    A^{T}U_{k+1} & = & V_{k}B_{k}^{T} + \alpha_{k+1}v_{k+1}e_{k+1}^{T}.
\end{eqnarray}
Also, the first equiality can be rewritten as
\begin{equation*}
    U^{T}_{k+1} A V_{k} = B_{k},
\end{equation*}
explaining the notion of left and right Lanczos vectors. Moreover, $ u_{k} \in
\mathcal{K}_{k}\left(AA^{T}, u_{1}\right)$, $v_{k} \in
\mathcal{K}_{k}\left(A^{T}A, v_{1}\right)$, where
\begin{eqnarray}
    \label{eqn:krlsubs1}
    \mathcal{K}_{k}\left(AA^{T}, u_{1}\right) = \left\{u_{1},
        AA^{T}u_{1}, \dots, \left(AA^{T}\right)^{k}u_{1}\right\}, \\
    \label{eqn:krlsubs2}
    \mathcal{K}_{k}\left(A^{T}A, v_{1}\right) = \left\{v_{1},
        A^{T}Av_{1}, \dots, \left(A^{T}A\right)^{k}v_{1}\right\},
\end{eqnarray}
and therefore $u_{1}, u_{2}, \dots, u_{k}$ and $v_{1}, v_{2}, \dots, v_{k}$ form
an orthogonal basis for these two Krylov subspaces.

\subsection{Lanczos Bidiagonalization in Inexact Arithmetics}
\label{sec:inexarith}
When the Lanczos bidiagonalization is carried out in finite precision
arithmetics the Lanczos recurrence becomes
\begin{eqnarray*}
    \alpha_{j} v_{j} & = & A^{T}u_{j} - \beta_{j}v_{j-1} + f_{j}, \\
    \beta_{j+1} u_{j+1} & = & Av_{j} -\alpha_{j} u_{j} + g_{j},
\end{eqnarray*}
where $f_{j} \in \mathbb{R}^{K}$, $g_{j} \in \mathbb{R}^{L}$ are error vectors
accounting for the rounding errors at the $j$-th step. Usually, the rounding
terms are small, thus after $k$ steps the equations \eqref{eqn:lbdiagm1},
\eqref{eqn:lbdiagm2} still hold to almost machine accuracy. In contrast, the
orthogonality among left and right Lanczos vectors is lost.

However, the Lanzos algorithm possesses a remarkable feature that the accuracy
of the converged Ritz values is not affected by the loss of orthogonality, but
while the largest singular values of $B_{k}$ are still accurate approximations
to the largest singular values of $A$, the spectrum of $B_{k}$ will in addition
contain false multiple copies of converged Ritz values (this happens even if the
corresponding true singular values of $A$ are isolated). Moreover, spurious
singular values (``ghosts'') periodically appear between the already converged
Ritz values (see, e.g. \cite{Larsen98} for illustration).

There are two different approaches with respect to what should be done to obtain a
robust method in finite arithmetics.

\subsubsection{Lanczos Bidiagonalization with no orthogonalization}
One approach is to apply the simple Lanczos process ``as-is'' and subsequently
use some criterion to detect and remove spurious singular values.

The advantage of this method is that it completely avoids any extra work
connected with the reorthogonalization, and the storage requirements are very
low since only few of the latest Lanczos vectors have to be remembered. The
disadvantage is that many iterations are wasted on simply generating multiple
copies of the large values \cite{Parlett80}. The number of extra iterations
required compared to that when executing the Lanczos process in exact
arithmetics can be very large: up to six times the original number as reported
in \cite{Parlett81}. Another disdvantage is that the criterion mentioned
above can be rather difficult to implement, since it depends on the correct
choice of different thresholds and knobs.

\subsubsection{Lanczos Bidiagonalization using reorthogonalization}
A different way to deal with loss of orthogonalizy among Lanczos vectors is to
apply some \emph{reorthogonalization} scheme. The simplest one is so-callled
\emph{full reorthogonalization} (FRO) where each new Lanczos vector $u_{j+1}$
(or $v_{j+1}$) is orthogonalized against all the Lanczos vectors generated so
far using, e.g., the modified Gram-Schmidt algorithm.

This approach is usually too expensive, since for the large problem the running
time needed for reorthogonalization quickly dominate the execution time unless
the necessary number of iterations $k$ is very small compared to the dimensions
of the problem (see section \ref{ssec:lancomplexity} for full analysis of the
complexity). The storage requirements of FRO might also be a limiting factor
since all the generated Lanczos vectors need to be saved.

A number of different schemes for reducing the work associated with keeping the
Lanczos vectors orthogonal were developed (see \cite{Baglama05} for extensive
review). The main idea is to estimate the \emph{level of orthogonality} among
the Lanczos vectors and perform the orthogonalization only when necessary.

For example, the so-called \emph{partial reorthogonalization scheme} (PRO)
described in \cite{Simon84a, Simon84b, Larsen98} uses the fact that the
level of orthogonality among the Lanczos vectors satisfies some recurrence
relation which can be derived from the recurrence used to generate the vectors
themselves.

\subsection{The Complexity of the Lanczos Bidiagonalization}
\label{ssec:lancomplexity}
The computational complexity of the Lanczos bidiagonalization is determined by
the two matrix-vector multiplications, thus the overall complexity of $k$
Lanczos iterations in exact arithmetics is $O(kLK)$ multiplications.

The computational costs increase when one deals with the lost of orthogonality
due to the inexact computations. For FRO scheme the additional $O(k^{2}(L+K))$
operations required for the orthogonalization quickly dominate the execution
time, unless the necessary number of iterations $k$ is very small compared to
the dimension of the problem. The storage requirements of FRO may also be a
limiting factor since all the generated Lanczos vectors need to be saved.

In general, there is no way to determine in advance how many steps will be
needed to provide the singular values of interest within a specified
accuracy. Usually the number of steps is determined by the distribution of the
singular numbers and the choice of the starting vector $p_{0}$. In the case when
the singular values form a cluster the convergence will not occur until $k$, the
number of iterations, gets very large. In such situation the method will also
suffer from increased storage requirements: it is easy to see that they are
$O(k(L+K))$.

These problems can be usually solved via limiting the dimension of the Krylov
subspaces \eqref{eqn:krlsubs1}, \eqref{eqn:krlsubs2} and using \emph{restarting
    schemes}: restarting the iterations after a number of steps by replacing the
starting vector with some ``improved'' starting vector (ideally, we would like
$p_{0}$ to be a linear combination of the right singular vectors of $A$
associated with the desired singular values). Thus, if we limit the dimension of
the Krylov subspaces to $d$ the storage complexity drops to $O(d(L+K))$ (surely
the number of desired singular triples should be greater than $d$). See
\cite{Baglama05} for the review of the restarting schemes proposed so far and
introduction to the thick-restarted Lanczos bidiagonalization algorithm. In the
papers \cite{Wu00, Yamazaki08} different strategies of choosing the best
dimension $d$ of Krylov subspaces are discussed.

If the matrix being decomposed is sparse or structured (for example, Hankel in
Basic SSA case or Hankel with Hankel blocks for 2D-SSA) then the computational
complexity of the bidiagonalization can be significantly reduced given that the
efficient matrix-vector multiplication routine is available.

\section{Efficient Matrix-Vector Multiplication for Hankel Matrices}
\label{sec:matmul}
In this section we will present efficient matrix-vector multiplication algorithm
for Hankel SSA trajectory matrix. By means of the fast Fourier transform the
cost of the multiplication drops from ordinary $O(LK)$ multiplications to
$O\left((L+K)\log(L+K)\right)$. The idea of using the FFT for computing a Hankel
matrix-vector product is nothing new: it was probably first described by
Bluestein in 1968 (see \cite{Swarztrauber91}). A later references are
\cite{OLeary81, Browne09}.

This approach reduces the computational complexity of the singular value
decomposition of Hankel matrices from $O\left(kLK + k^{2}(L+K)\right)$ to
$O\left(k(L+K)\log(L+K) + k^{2}(L+K)\right)$. The complexity of the worst case
when $K \sim L \sim N/2$ drops significantly from $O\left(kN^{2} +
    k^{2}N\right)$ to $O\left(kN\log{N} + k^{2}N\right)$\footnote{Compare with
    $O\left(N^{3}\right)$ for full SVD via Golub-Reinsh method.}\!\!.

Also the discribed algorithms provide an efficient way to store the Hankel SSA
trajectory matrix reducing the storage requirements from $O(LK)$ to $O(L+K)$.

\subsection{Matrix Representation of Discrete Fourier Transform and Circulant
    Matrices}
\label{ssec:fftcirc}
The 1-d discrete Fourier transform (DFT) of a (complex) vector
$\left(f_{k}\right)_{k=0}^{N-1}$ is defined by:
\begin{equation*}
    \hat{f}_{l} = \sum_{k=0}^{N}{e^{-2 \pi i k l / N}f_{k}}, \quad k =
    0,\dots, N-1.
\end{equation*}
Denote by $\omega$ the primitive $N$-root of unity, $\omega = e^{-2 \pi i / N}$.
Introduce the DFT matrix $\mathbb{F}_{N}$:
\begin{equation*}
    \mathbb{F}_{N} =
    \begin{pmatrix}
        1      & 1           & \cdots & 1                 \\
        1      & \omega^{1}  & \cdots & \omega^{N-1}       \\
        1      & \omega^{2}  & \cdots & \omega^{2(N-1)}     \\
        \vdots & \vdots      & \vdots & \vdots            \\
        1      & \omega^{N-1} & \cdots & \omega^{(N-1)(N-1)}
    \end{pmatrix}.
\end{equation*}
Then the 1-d DFT can be written in matrix form: $\hat{f} = \mathbb{F}_{N}f$. The
inverse of the DFT matrix is given by
\begin{equation*}
    \mathbb{F}_{N}^{-1} = \frac{1}{N}\mathbb{F}_{N}^{*},
\end{equation*}
where $\mathbb{F}_{N}^{*}$ is the adjoint (conjugate transpose) of DFT matrix
$\mathbb{F}_{N}$. Thus, the inverse 1-d discrete Fourier transform (IDFT) is
given by
\begin{equation*}
    f_{k} = \frac{1}{N}\sum_{l=0}^{N-1}{e^{2 \pi i k l / N}\hat{f}_{l}},
    \quad k = 0,\dots,N-1.
\end{equation*}

The fast Fourier transform is an efficient algorithm to compute the DFT (and
IDFT) in $O(N\log{N})$ complex multiplications instead of $N^{2}$ complex
multiplications in direct implementation of the DFT \cite{Loan92, FFTW05}.

\begin{definition}
    An $N \times N$ matrix $C$ of the form
    \begin{equation*}
        C =
        \begin{pmatrix}
            c_{1}   & c_{N}   & c_{N-1} & \cdots & c_{3}   & c_{2} \\
            c_{2}   & c_{1}   & c_{N}   & \cdots & c_{4}   & c_{3} \\
            \vdots  & \vdots & \vdots &        & \vdots  & \vdots \\
            c_{N-1} & c_{N-2} & c_{N-3} & \cdots & c_{1}   & c_{N} \\
            c_{N}   & c_{N-1} & c_{N-2} & \cdots & c_{2}   & c_{1}
        \end{pmatrix}
    \end{equation*}
    is called a \bf{circulant matrix}.
\end{definition}

A circulant matrix is fully specified by its first column $c =
\left(c_{1},c_{2},\dots,c_{N}\right)^{T}$ The remaining columns of $C$ are each
cyclic permutations of the vector $c$ with offset equal to the column index. The
last row of $C$ is the vector $c$ in reverse order, and the remaining rows are
each cyclic permutations of the last row.

The eigenvectors of a circulant matrix of a given size are the columns of the
DFT matrix of the same size. Thus a circulant matrix is diagonalized by the DFT
matrix:
\begin{equation*}
    C = \mathbb{F}_{N}^{-1} \diag(\mathbb{F}_{N}c) \mathbb{F}_{N},
\end{equation*}
so the eigenvalues of $C$ are given by the product $\mathbb{F}_{N}c$.

This factorization can be used to perform efficient matrix-vector
multiplication. Let $v \in \mathbb{R}^{N}$ and $C$ be an $N \times N$ circulant
matrix. Then
\begin{equation*}
    Cv = \mathbb{F}_{N}^{-1} \diag(\mathbb{F}_{N}c) \mathbb{F}_{N} v =
    \mathbb{F}_{N}^{-1} \left(\mathbb{F}_{N}c \odot \mathbb{F}_{N}v\right),
\end{equation*}
where $\odot$ denotes the element-wise vector multiplication. This can be
computed efficiently using the FFT by first computing the two DFTs,
$\mathbb{F}_{N}c$ and $\mathbb{F}_{N}v$, and then computing the IDFT
$\mathbb{F}_{N}^{-1} \left(\mathbb{F}_{N}c \odot \mathbb{F}_{N}v\right)$. If the
matrix-vector multiplication is performed repeatedly with the same circulant
matrix and different vectors, then surely the DFT $\mathbb{F}_{N}c$ needs to be
computed only once.

In this way the overall computational complexity of matrix-vector product for
the circulant matrices drops from $O(N^{2})$ to $O(N\log{N})$.

\subsection{Toeplitz and Hankel Matrices}
\label{sec:hmatmul}
\begin{definition}
    A $L \times K$ matrix of the form
    \begin{equation*}
        T =
        \begin{pmatrix}
            t_{K}    & t_{K-1}   & \cdots  & t_{2}  & t_{1}  \\
            t_{K+1}   & t_{K}    & \cdots  & t_{3}  & t_{2}  \\
            \vdots   & \vdots   &         & \vdots & \vdots \\
            t_{K+L-1} & t_{K+L-2} & \cdots  & t_{L+1} & t_{L}
        \end{pmatrix}
    \end{equation*}
    is called a (non-symmetric) \textbf{Toeplitz matrix}.
\end{definition}
A Toeplitz matrix is completely determined by its last column and last row, and
thus depends on $K+L-1$ parameters. The entries of $T$ are constant down the diagonals
parallel to the main diagonal.

Given an algorithm for the fast matrix-vector product for circulant matrices, it
is easy to see the algorithm for a Toeplitz matrix, since a Toeplitz matrix can
be embedded into a circulant matrix $C_{K+L-1}$ of size $K+L-1$ with the first
column equals to
$\left(t_{K},t_{K+1},\dots,t_{K+L-1},t_{1},t_{2},\dots,t_{K-1}\right)^{T}$:
\begin{equation}
    \label{eqn:embcirc}
    C_{K+L-1} =
    \begin{pmatrix}
        t_{K}      & \cdots & t_{1} & t_{K+L-1} & \cdots & t_{K+1} \\
        t_{K+1}    & \cdots & t_{2} & t_{1} &  \cdots & t_{K+2} \\
        t_{K+2}    & \cdots & t_{3} & t_{2} &  \cdots & t_{K+3} \\
        \vdots    &        & \vdots & \vdots &  & \vdots \\
        t_{K+L-1}  & \cdots & t_{L} & t_{L-1} & \cdots & t_{1} \\
        t_{1}     & \cdots & t_{L+1} & t_{L} & \cdots & t_{2} \\
        t_{2}     & \cdots & t_{L+2} & t_{L+1} & \cdots & t_{3} \\
        \vdots    &        & \vdots & \vdots &  & \vdots \\
        t_{K-1}   & \cdots & t_{K+L-1} & t_{K+L-2} & \cdots & t_{K} \\
    \end{pmatrix}
    ,
\end{equation}
where the leading $L \times K$ principal submatrix is $T$. This technique of
embedding a Toeplitz matrix in a larger circulant matrix to achieve fast
computation is widely used in preconditioning methods \cite{Chan93}.

Using this embeddings, the Toeplitz matrix-vector product $Tv$ can be
represented as follows:
\begin{equation}
    \label{eqn:tmatmul}
    C_{K+L-1}
    \begin{pmatrix}
        v \\
        \mathbf{0}_{L-1}
    \end{pmatrix} =
    \begin{pmatrix}
        Tv \\
        *
    \end{pmatrix},
\end{equation}
and can be computed efficiently in $O\left((K+L)\log(K+L)\right)$ time and
$O(K+L)$ memory using the FFT as described previously.

Recall that an $L \times K$ matrix of form
\begin{equation*}
    H =
    \begin{pmatrix}
        h_{1}    & h_{2}   & \cdots  & h_{K-1}  & h_{K}  \\
        h_{2}    & h_{3}   & \cdots  & h_{K}    & h_{K+1}  \\
        \vdots  & \vdots  &         & \vdots   & \vdots \\
        h_{L}    & h_{L+1} & \cdots  & h_{K+L-2} & h_{K+L-1}
    \end{pmatrix}
\end{equation*}
is called \textbf{Hankel matrix}.

A Hankel matrix is completely determined by its first column and last row, and
thus depends on $K+L-1$ parameters. The entries of $H$ are constant along the
antidiagonals. One can easily convert Hankel matrix into a Toeplitz one by
reversing its columns. Indeed, define
\begin{equation*}
    P =
    \begin{pmatrix}
        0 & 0 & \cdots & 0 & 1 \\
        0 & 0 & \cdots & 1 & 0 \\
        \vdots & \vdots & & \vdots & \vdots \\
        0 & 1 & \cdots & 0 & 0 \\
        1 & 0 & \cdots & 0 & 0
    \end{pmatrix}
\end{equation*}
the backward identity permutation matrix. Then $H P$ is a Toeplitz matrix for
any Hankel matrix $H$, and $T P$ is a Hankel matrix for any Toeplitz matrix
$T$. Note that $P = P^{T} = P^{-1}$ as well. Now for the product $Hv$ of Hankel
matrix $H$ and vector $v$ we have
\begin{equation*}
    Hv = HPPv = \left(HP\right)Pv = Tw,
\end{equation*}
where $T$ is Toeplitz matrix and vector $w$ is obtained as vector $v$ with
entries in reversed order. The product $Tw$ can be computed using circulant
embedding procedure as described in \eqref{eqn:tmatmul}.

\subsection{Hankel SSA Trajectory Matrices}
\label{ssec:hssamatmul}

Now we are ready to exploit the connection between the time series
$F=\left(f_{j}\right)_{j=1}^{N}$ under decomposition and corresponding
trajectory Hankel matrix to derive the matrix-vector multiplication algorithm in
terms of the series itself. In this way we will effectively skip the embedding
step of the SSA algorithm and reduce the computational and storage complexity.

The entries of the trajectory matrix of the series are $h_{j} =
f_{j},\,j=1,\dots,N$. The entries of Toeplitz matrix $T = HP$ are $t_{j} = h_{j}
= f_{j}$, $j=1,\dots,N$. The corresponding first column of the embedding
circulant matrix \eqref{eqn:embcirc} is
\begin{equation*}
    c_{j} =
    \begin{cases}
        t_{N-L+j} = f_{N-L+j}, & \quad 1 \leqslant j \leqslant L; \\
        t_{j-L+1} = f_{j-L+1} & \quad L < j \leqslant N,
    \end{cases}
\end{equation*}
and we end with the following algorithm for matrix-vector multiplication for
trajectory matrix:

\begin{algorithm}[H]
    \label{alg:hmatmul}
    \caption{Matrix-Vector Multiplication for SSA Trajectory Matrix}
    \dontprintsemicolon
    \KwIn{Time series $F=\left(f_{j}\right)_{j=1}^{N}$, window length $L$,
        vector $v$ of length $N-L+1$.}
    \KwOut{$p=Xv$, where $X$ is a trajectory matrix for $F$ given window length
        $L$.}
    \BlankLine
    $c \longleftarrow
    \left(f_{N-L+1},\dots,f_{N},f_{1},\dots,f_{N-L}\right)^{T}$ \nllabel{hmatmuls1} \;
    $\hat{c} \longleftarrow FFT_{N}(c)$ \nllabel{hmatmuls2} \;
    $w \longleftarrow \left(v_{N-L+1},\dots,v_{1},0,\dots,0\right)^{T}$ \;
    $\hat{w} \longleftarrow FFT_{N}(w)$ \;
    $p' \longleftarrow IFFT_{N}(\hat{w}\odot\hat{c})$ \;
    $p \longleftarrow \left(p'_{1},\dots,p'_{L}\right)^{T}$
\end{algorithm}
where $(\hat{w}\odot\hat{c})$ denotes element-wise vector multiplication.

If the matrix-vector multiplication $Xv$ is performed repeatedly with the same
matrix $X$ and different vectors $v$ then steps \ref{hmatmuls1}, \ref{hmatmuls2}
should be performed only once at the beginning and the resulting vector
$\hat{c}$ should be reused later on.

The overall computational complexity of the matrix-vector multiplication is
$O(N\log{N})$. Memory space of size $O(N)$ is required to store precomputed
vector $\hat{c}$.

Note that the matrix-vector multiplication of the transposed trajectory matrix
$X^{T}$ can be performed using the same vector $\hat{c}$. Indeed, the
bottom-right $K \times L$ submatrix of circulant matrix \eqref{eqn:embcirc}
contains the Toeplitz matrix $X^{T}P$ and we can easily modify algorithm
\ref{alg:hmatmul} as follows:

\begin{algorithm}[H]
    \caption{Matrix-Vector Multiplication for transpose of SSA Trajectory Matrix}
    \dontprintsemicolon
    \KwIn{Time series $F=\left(f_{j}\right)_{j=1}^{N}$, window length $L$,
        vector $v$ of length $L$.}
    \KwOut{$p=X^{T}v$, where $X$ is a trajectory matrix for $F$ given window
        length $L$.}
    \BlankLine
    $c \longleftarrow
    \left(f_{N-L+1},\dots,f_{N},f_{1},\dots,f_{N-L}\right)^{T}$ \;
    $\hat{c} \longleftarrow FFT_{N}(c)$ \ \;
    $w \longleftarrow \left(0,\dots,0,v_{L},\dots,v_{1}\right)^{T}$ \;
    $\hat{w} \longleftarrow FFT_{N}(w)$ \;
    $p' \longleftarrow IFFT_{N}(\hat{w}\odot\hat{c})$ \;
    $p \longleftarrow \left(p'_{L-1},\dots,p'_{N}\right)^{T}$
\end{algorithm}

\section{Efficient Rank 1 Hankelization}
\label{sec:diagavg}
Let us recall the \emph{diagonal averaging} procedure as described in
\ref{ssec:hankelization} which transforms an arbitrary $L \times K$ matrix $Y$
into a Hankel one and therefore into a series $g_{k}$, $k=1, \dots,N$.
\begin{equation}
    \label{eqn:hankelization}
    g_{k} =
    \begin{cases}
        \frac{1}{k}\sum_{j=1}^{k}{y_{j,k-j+1}}, & \quad 1 \leqslant k \leqslant L, \\
        \frac{1}{L}\sum_{j=1}^{L}{y_{j,k-j+1}}, & \quad L < k < K, \\
        \frac{1}{N-k+1}\sum_{j=k-K+1}^{L}{y_{j,k-j+1}}, & \quad K \leqslant k
        \leqslant N.
    \end{cases}
\end{equation}

Without loss of generality we will consider only \emph{rank 1 hankelization}
when matrix $Y$ is elementary and can be represented as $Y =
\sigma uv^{T}$ with vectors $u$ and $v$ of length $L$ and $K$
correspondingly. Then equation \eqref{eqn:hankelization} becomes
\begin{equation}
    \label{eqn:hankelization1}
    g_{k} =
    \begin{cases}
        \frac{\sigma}{k}\sum_{j=1}^{k}{u_{j}v_{k-j+1}}, & \quad 1 \leqslant k \leqslant L, \\
        \frac{\sigma}{L}\sum_{j=1}^{L}{u_{j}v_{k-j+1}}, & \quad L < k < K, \\
        \frac{\sigma}{N-k+1}\sum_{j=k-K+1}^{L}{u_{j}v_{k-j+1}}, & \quad K \leqslant k
        \leqslant N.
    \end{cases}
\end{equation}
This gives a hint how the diagonal averaging can be efficiently computed.

Indeed, consider the infinite series $u'_{n}$, $n \in \mathbb{Z}$ such that
$u'_{n} = u_{n}$ when $1 \leqslant n \leqslant L$ and $0$ otherwise. The series
$v'_{n}$ are defined in the same way, so $v'_{n} = v_{n}$ when $1 \leqslant n
\leqslant K$ and $0$ otherwise. The linear convolution $u'_{n} * v'_{n}$ of
$u'_{n}$ and $v'_{n}$ can be written as follows:
\begin{multline}
    \label{eqn:conv}
    \left(u'_{n} * v'_{n}\right)_{k} =
    \sum_{j=-\infty}^{+\infty}{u'_{j}v'_{k-j+1}} =
    \sum_{j=1}^{L}{u'_{j}v'_{k-j+1}} = \\ =
    \begin{cases}
        \sum_{j=1}^{k}{u_{j}v_{k-j+1}}, & \quad 1 \leqslant k \leqslant L, \\
        \sum_{j=1}^{L}{u_{j}v_{k-j+1}}, & \quad L < k < K, \\
        \sum_{j=k-K+1}^{L}{u_{j}v_{k-j+1}}, & \quad K \leqslant k
        \leqslant N.
    \end{cases}
\end{multline}
Comparing the equations \eqref{eqn:hankelization} and \eqref{eqn:conv} we deduce
that
\begin{equation*}
    g_{k} = c_{k} \left(u'_{n} * v'_{n}\right)_{k},
\end{equation*}
where the constants $c_{k}$ are known in advance.

The linear convolution $u * v$ can be defined in terms of the \emph{circular
    convolution}, as follows. Pad the vectors $u$ and $v$ with the zeroes up to
length $L + K -1$. Then the linear convolution of the original vectors $u$ and
$v$ is the same as the circular convolution of the extended series of length
$L+K-1$, namely
\begin{equation*}
    \left(u \ast v\right)_{k} = \sum_{j=1}^{L+K-1}{u_{j}v_{k-j+1}},
\end{equation*}
where $k-j+1$ is evaluated $\mathrm{mod}\,(L+K-1)$. The resulting circular
convolution can be calculated effeciently via the FFT \cite{Briggs95}:
\begin{equation*}
    \left(u \ast v\right)_{k} =
    IFFT_{N}\left(FFT_{N}(u') \odot FFT_{N}(v')\right).
\end{equation*}
Here $N = L+K-1$, and $u'$, $v'$ denote the zero-padded vectors $u$ and $v$
correspondingly.

And we end with the following algorithm for the hankelization via
convolution:

\begin{algorithm}[H]
    \caption{Rank 1 Hankelization via Linear Convolution}
    \dontprintsemicolon
    \KwIn{Vector $u$ of length $L$, vector $v$ of length $K$, singular value $\sigma$.}
    \KwOut{Time series $G=\left(g_{j}\right)_{j=1}^{N}$ corresponding to the
        matrix $Y = \sigma uv^{T}$ after hankelization.}
    \BlankLine
    $u' \longleftarrow
    \left(u_{1},\dots,u_{L},0,\dots,0\right)^{T} \in \mathbb{R}^{N}$ \;
    $\hat{u} \longleftarrow FFT_{N}(u')$ \;
    $v' \longleftarrow
    \left(v_{1},\dots,v_{K},0,\dots,0\right)^{T} \in \mathbb{R}^{N}$ \;
    $\hat{v} \longleftarrow FFT_{N}(v')$ \;
    $g' \longleftarrow IFFT_{N}(\hat{v}\odot\hat{u})$ \;
    $w \longleftarrow \left(1,\dots,L,L,\dots,L,L,\dots,1\right) \in
    \mathbb{R}^{N}$ \;
    $g \longleftarrow \sigma \left(w \odot g'\right)$
\end{algorithm}
The computational complexity of this algorithm is $O(N\log{N})$ multiplications
versus $O(LK)$ for na\"\i ve direct approach. Basically this means that the
worst case complexity of the diagonal averaging drops from $O(N^{2})$ to
$O(N\log{N})$.

\section{Implementation Comparison}
\label{sec:impl}
The algorithms presented in this paper have much better asymptotic complexity
than standard ones, but obviously might have much bigger overhead and thus would
not be suitable for the real-world problems unless series length and/or window
length becomes really big.

The mentioned algorithms were implemented by the means of the
\verb@Rssa@\footnote{The package will be submitted to CRAN soon, yet it can be
    obtained from author.} package for the \verb@R@ system for statistical
computing \cite{R09}. All the comparions were carried out on the 4 core AMD
Opteron 2210 workstation running Linux. Where possible we tend to use
state-of-the-art implementations of various computational algorithms
(e.g. highly optimized \verb@ATLAS@ \cite{Whaley04} and \verb@ACML@
implementations of \verb@BLAS@ and \verb@LAPACK@ routines). We used \verb@R@
version 2.8.0 throughout the tests.

We will compare the performance of such key components of the fast SSA
implementation as rank 1 diagonal averaging and the Hankel matrix-vector
multiplication. The performance of the methods which use the ``full'' SSA,
namely the series bootstrap confidence interval construction and series
structure change detection, will be considered as well.

\subsection{Fast Hankel Matrix-Vector Product}
\label{ssec:hmatmulcmp}
Fast Hankel matrix-vector multiplication is the key component for the Lanczos
bidiagonalization. The algorithm outlined in the section \ref{ssec:hssamatmul}
was implemented in two different ways: one using the core R FFT implementation
and another one using FFT from FFTW library \cite{FFTW05}. We selected window
size equal to the half of the series length since this corresponds to the worst
case both in terms of computational and storage complexity. All initialization
times (for example, Hankel matrix construction or circulant precomputation
times) were excluded from the timings since they are performed only once in the
beginning. So, we compare the performance of the generic \verb@DGEMV@
\cite{Lapack99} matrix-vector multication routine versus the performance of the
special routine exploiting the Hankel structure of the matrix.

\begin{figure}[htbp]
    \centering
    \includegraphics[width=\columnwidth]{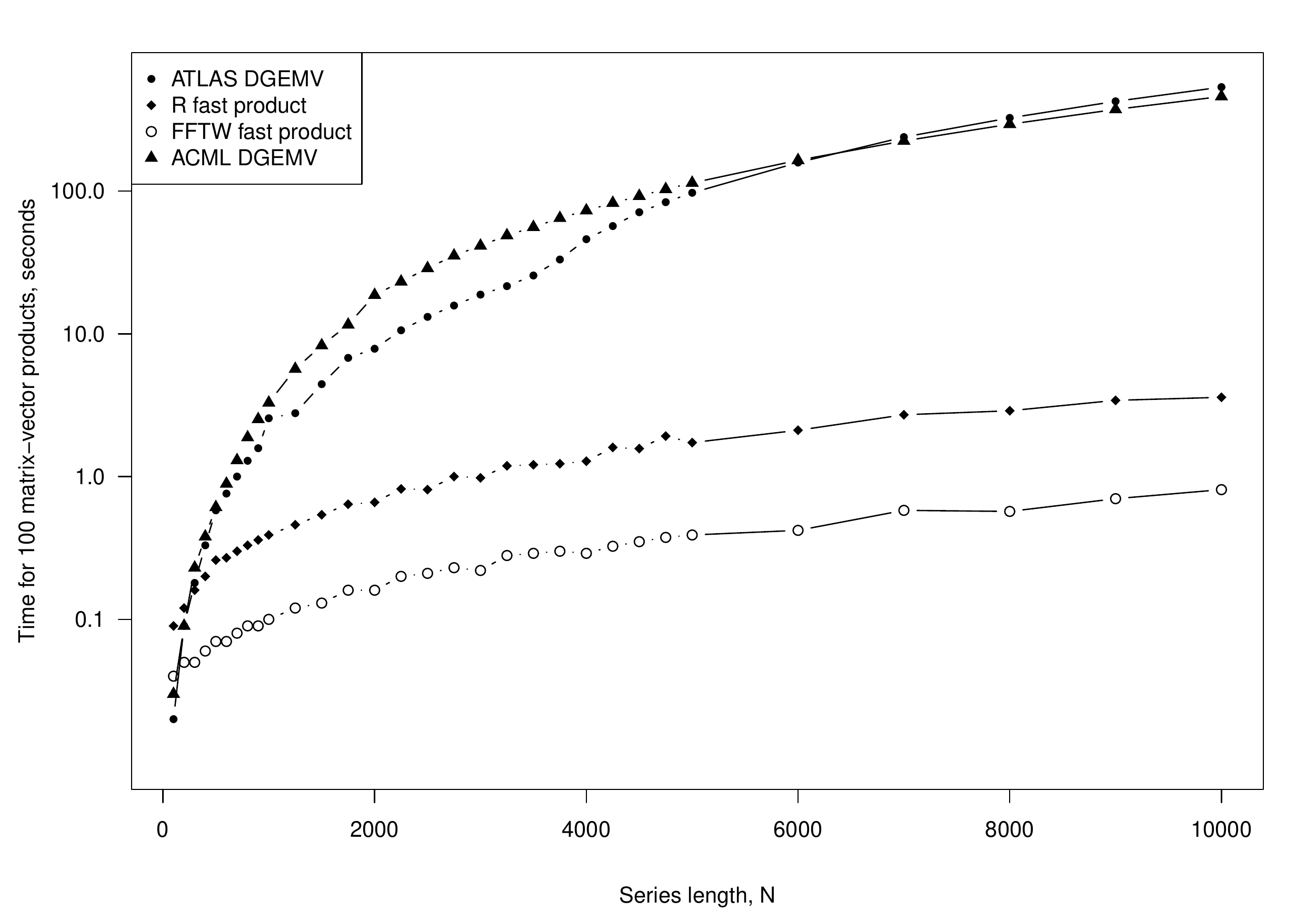}
    \caption{Hankel Matrix-Vector Multiplication Comparison}
    \label{fig:hmatmul}
\end{figure}

The running times for 100 matrix-vector products with different series lengths
are presented on the figure \ref{fig:hmatmul}. Note that we had to use
logarithmic scale for the y-axis in order to outline the difference between the
generic routines and our implementations properly. One can easily see that all
the routines indeed possess the expected asymptotic complexity behaviour. Also
special routines are almost always faster than generic ones, thus we recommend
to use them for all window sizes (the overhead for small series lengths can be
neglected).

\subsection{Rank 1 Hankelization}
\label{ssec:hankelizationcmp}
The data access during straightforward implementation of rank 1 hankelization is
not so regular, thus pure R implementation turned out to be slow and we haven't
considered it at all. In fact, two implementations are under comparison:
straightforward C implementation and FFT-based one as described in section
\ref{sec:diagavg}. For the latter FFTW library was used. The window size was
equal to the half of the series length to outline the worst possible case.

\begin{figure}[htbp]
    \centering
    \includegraphics[width=\columnwidth]{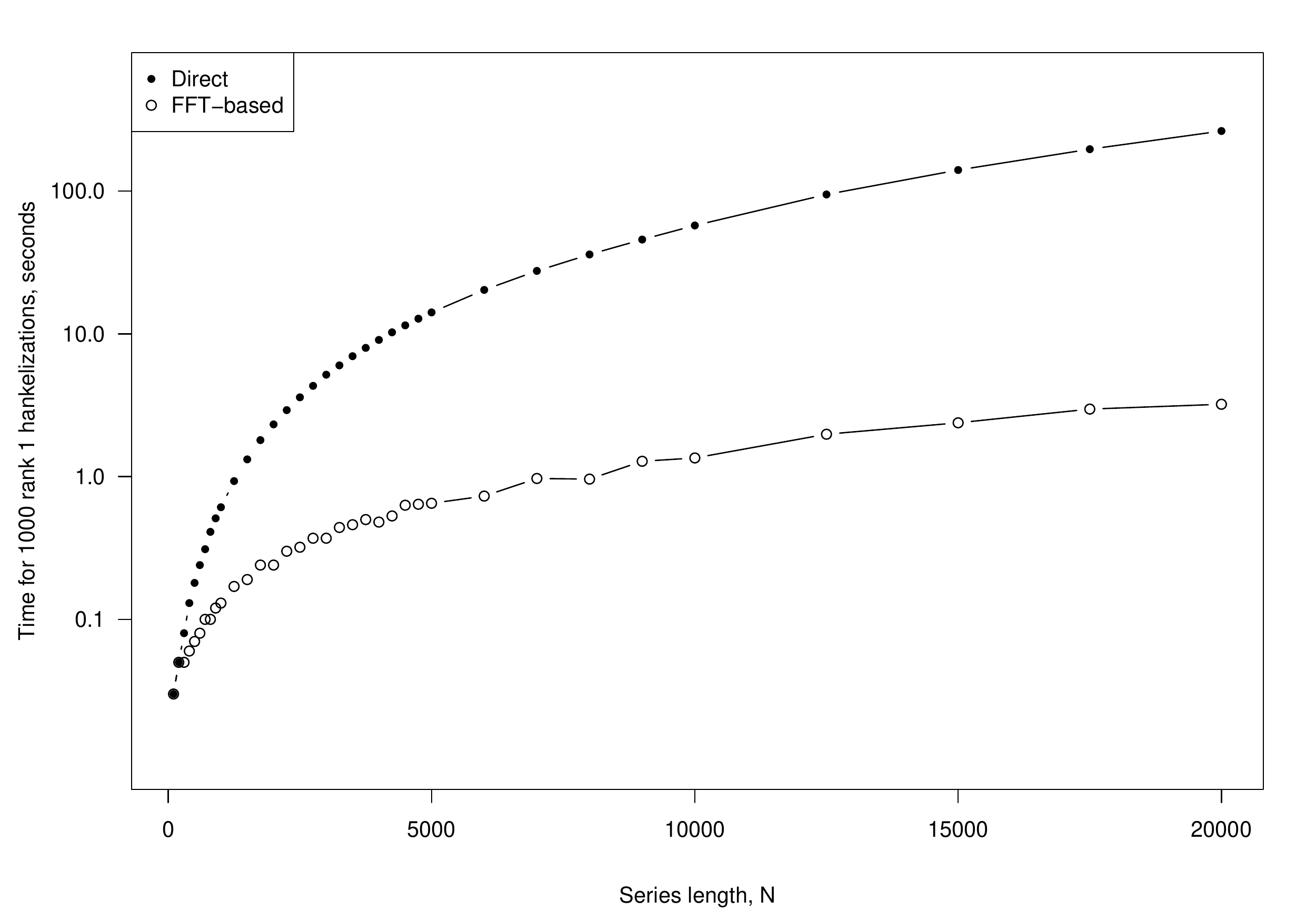}
    \caption{Rank 1 Hankelization Implementation Comparison}
    \label{fig:diagavg}
\end{figure}

The results are shown on the figure \ref{fig:diagavg}. As before, logarithmic
scale was used for the y-axis. From this fugure one can see that the
computational complexity of the direct implementation of hankelization quickly
dominates the overhead of the more complex FFT-based one and thus the latter can
be readily used for all non-trivial series lengths.

\subsection{Bootstrap Confidence Intervals Construction}
For the comparison of the implementations of the whole SSA algorithm the problem
of bootstrap confidence intervals construction was selected. It can be described
as follows: consider we have the series $F_{N}$ of finite rank $d$. Form the
series $F'_{N} = F_{N} + \sigma \varepsilon_{N}$, where $\varepsilon_{i}$
denotes the uncorrelated white noise (thus the series $F'_{N}$ are of full
rank). Fix the window size $L < N$ and let $G_{N}$ denote the series
reconstructed from first $d$ eigentriples of series $F'_{N}$. This way $G_{N}$
is a series of rank $d$ as well.

Since the original series $F_{N}$ are considered known, we can study
large-sample properties of the estimate $G_{N}$: bias, variance and mean square
error. Having the variance estimate one can, for example, construct bootstrap
confidence intervals for the mean of the reconstructed series $G_{N}$.

Such simulations are usually quite time-consuming since we need to perform few
dozen reconstructions for the given length $N$ in order to estimate the
variance.

For the simulation experiments we consider the series $F_{N} =
(f_{1},\dots,f_{N})$ of rank $5$ with
\begin{equation*}
    f_{n} = 10e^{-5n/N} +
    \sin\left(\frac{2\pi}{13}\frac{n}{N}\right) +
    2.5 \sin\left(\frac{2\pi}{37}\frac{n}{N}\right)
\end{equation*}
and $\sigma = 5$. The series $F_{1000}$ (white) and $F'_{1000}$ (black) are
shown on the figure \ref{fig:sinnoise}.

\begin{figure}[htbp]
    \centering
    \includegraphics[width=\columnwidth]{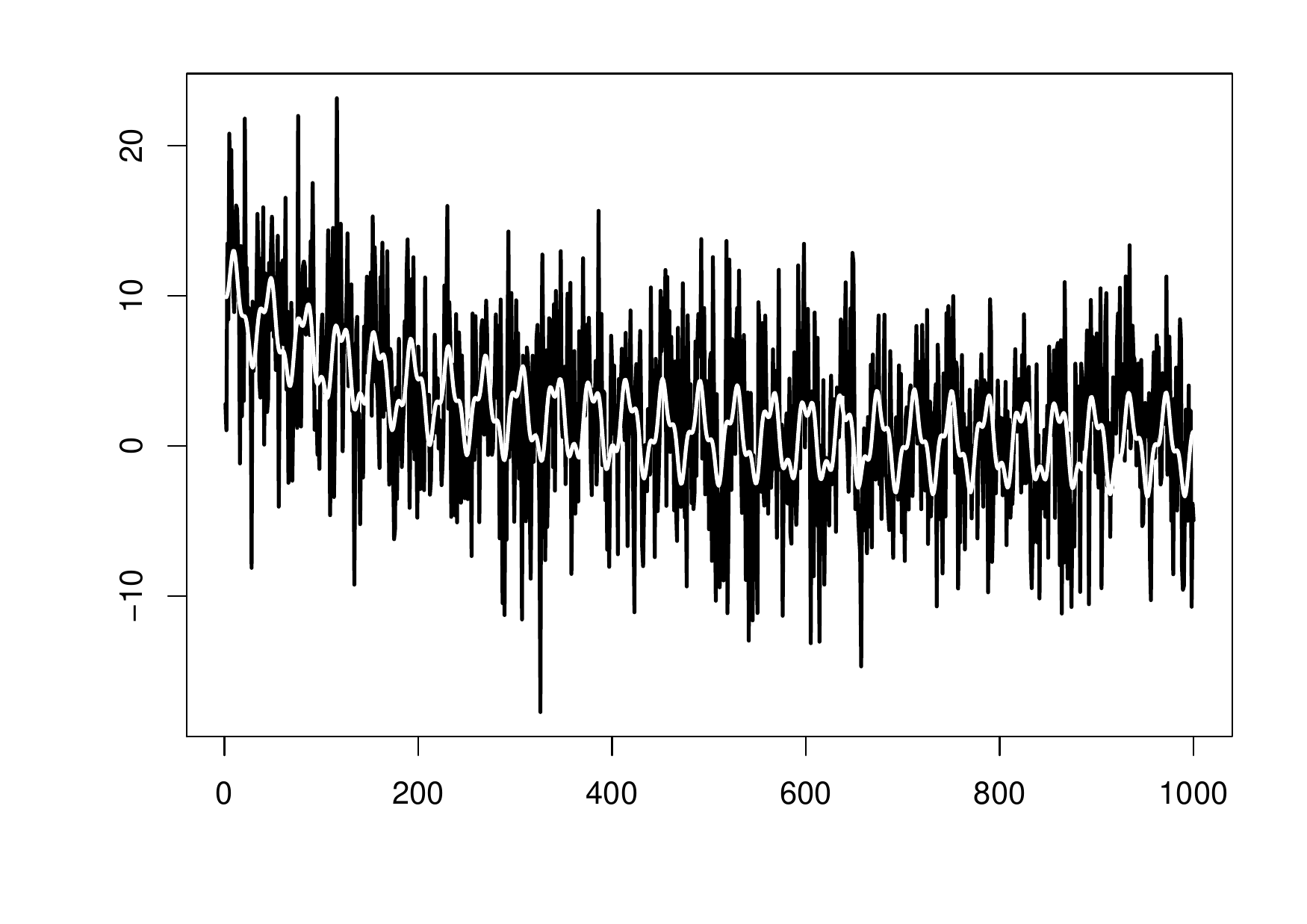}
    \caption{The Series $F_{1000}$ and $F'_{1000}$.}
    \label{fig:sinnoise}
\end{figure}

We compare $3$ different SSA implementations: one uses the full SVD via
\verb@DGESDD@ routine from \verb@ATLAS@ library, another one uses the full SVD
via the same routine from \verb@ACML@ library and third uses fast Hankel
truncated SVD and FFT-based rank-1 hankelization. Again, \verb@FFTW@ library was
used to perform the FFT. Window size for SSA was always fixed at half of the
series length.

\begin{figure}[htbp]
    \centering
    \includegraphics[width=\columnwidth]{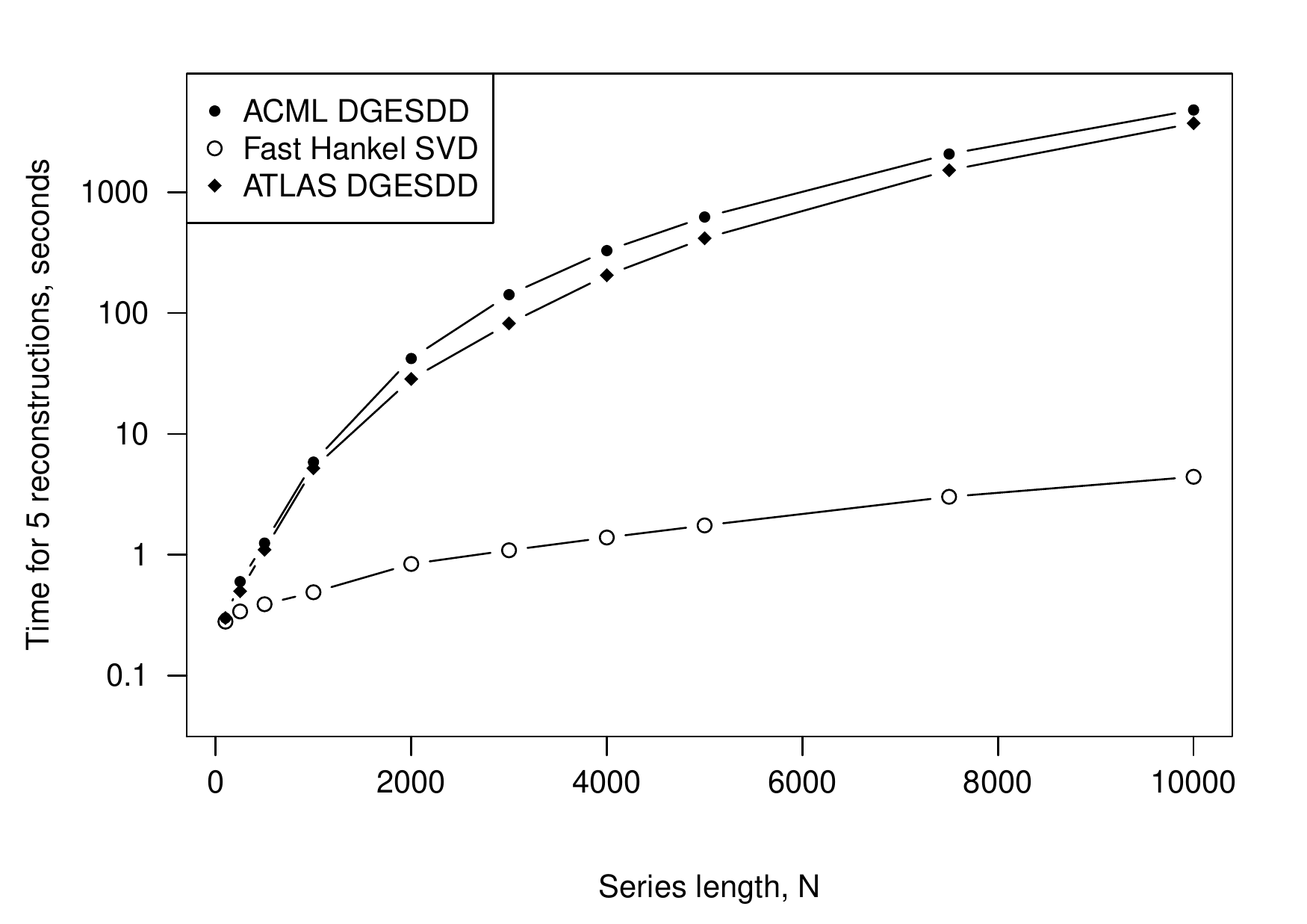}
    \caption{Reconstruction Time Comparison.}
    \label{fig:sinnoise}
\end{figure}

The running times for the rank $5$ reconstruction are presented on the figure
\ref{fig:sinnoise}. Note that logarithmic scale was used for the y-axis. From it
we can easily see that the running times for the SSA implementation via fast
Hankel SVD and FFT-based rank-1 hankelization are dramatically smaller compared
to the running times of other implementations for any non-trivial series
length.

\subsection{Detection of Structural Changes}
The SSA can be used to detect the structural changes in the series. The main
instrument here is so-called \emph{heterogeneity matrix} (H-matrix). We will
briefly describe the algorithm for construction of such matrix, discuss the
computation complexity of the contruction and compare the performances of
different implementations.

Exhaustive exposition of the detection of structural changes by the means of SSA
can be found in \cite{Golyandina01}.

Consider two time series $F^{(1)}$ and $F^{(2)}$. Let $N_{1}$ and $N_{2}$ denote
their lengths respectively. Take integer $L < \min\left(N_{1}-1,
    N_{2}\right)$. Let $U_{j}^{(1)}$, $j=1,\dots,L$ denote the eigenvectors of
the SVD of the trajectory matrix of the series $F^{(1)}$. Fix $I$ to be a
subset of $\left\{1,\dots,L\right\}$ and denote $\mathcal{L}^{(1)} =
\vspan\left(U_{i},\, i \in I\right)$. Denote by
$X^{(2)}_{1},\dots,X^{(2)}_{K_{2}}$ ($K_{2} = N_{2} - L + 1$) the $L$-lagged
vectors of the series $F^{(2)}$.

Now we are ready to introduce the \emph{heterogeneity index} \cite{Golyandina01}
which measures the descrepancy between the series $F^{(2)}$ and the structure of
the series $F^{(1)}$ (this structure is described by the subspace
$\mathcal{L}^{(1)}$):
\begin{equation}
    \label{eqn:g}
    g(F^{(1)},F^{(2)}) =
    \frac{\sum_{j=1}^{K_{2}}{\dist^{2}\left(X^{(2)}_{j}, \mathcal{L}^{(1)}\right)}}
    {\sum_{j=1}^{K_{2}}{\left\|X^{(2)}_{j}\right\|^{2}}},
\end{equation}
where $\dist(X,\mathcal{L})$ denotes the Euclidian distance between the vector
$X$ and subspace $\mathcal{L}$. One can easily see that $0 \leq g \leq 1$.

Note that since the vectors $U_{i}$ form the orthonormal basis of the subspace
$\mathcal{L}^{(1)}$ equation \eqref{eqn:g} can be rewritten as
\begin{equation*}
    g(F^{(1)},F^{(2)}) =
    1-
    \frac{\sum_{j=1}^{K_{2}}{\sum_{i \in I}}{\left(U_{i}^{T}X^{(2)}_{j}\right)^{2}}}
    {\sum_{j=1}^{K_{2}}{\left\|X^{(2)}_{j}\right\|^{2}}}.
\end{equation*}

Having this descrepancy measure for two arbitrary series in the hand one can
obviously construct the method of detection of structural changes in the single
time series. Indeed, it is sufficient to calculate the heterogeneity index $g$
for different pairs of subseries of the series $F$ and study the obtained
results.

The \emph{heterogeneity matrix} (H-matrix) \cite{Golyandina01} provides a
consistent view over the structural descrepancy between different parts of the
series. Denote by $F_{i,j}$ the subseries of $F$ of the form: $F_{i,j} =
\left(f_{i},\dots,f_{j}\right)$. Fix two integers $B > L$ and $T \geq L$. These
integers will denote the lengths of \emph{base} and \emph{test} subseries
correspondingly. Introduce the H-matrix $G_{B,T}$ with the elements $g_{ij}$ as
follows:
\begin{equation*}
    g_{ij} = g(F_{i,i+B}, F_{j,j+T}),
\end{equation*}
for $i=1,\dots,N-B+1$ and $j=1,\dots,N-T+1$. In simple words we split the series
$F$ into subseries of lengths $B$ and $T$ and calculate the heteregoneity index
between all possible pairs of the subseries.

The computation complexity of the calculation of the H-matrix is formed by the
complexity of the $SVDs$ for series $F_{i,i+B}$ and complexity of calculation of
all heteregeneity indices $g_{ij}$ as defined by equation \eqref{eqn:g}.

The worst-case complexity of the SVD for series $F_{i,i+B}$ corresponds to the
case when $L \sim B/2$. Then the single decomposition costs $O(B^{3})$ for full
SVD via Golub-Reinsh method and $O(kB\log{B}+k^{2}B)$ via fast Hankel SVD as
presented in sections \ref{sec:svd} and \ref{sec:matmul}. Here $k$ denotes the
number of elements in the index set $I$. Since we need to make $N-B+1$
decompositions a total we end with the following complexities of the worst case
(when $B \sim N/2$): $O(N^{4})$ for full SVD and $O(kN^{2}\log{N}+k^{2}N^{2})$
for fast Hankel SVD.

Having all the desired eigenvectors $U_{i}$ for all the subseries $F_{i,i+B}$ in
hand one can calculate the H-matrix $G_{B,T}$ in $O(kL(T-L+1)(N-T+1)+(N-T+1)L)$
multiplications. The worst case corresponds to $L \sim T/2$ and $T \sim N/2$
yielding the $O(kN^{3})$ complexity for this step. We should note that this
complexity can be further reduced with the help of fast Hankel matrix-vector
multiplication, but for the sake of simplicity we won't present this result
here.

Summing the complexities we end with $O(N^{4})$ multiplications for H-matrix
computation via full SVD and $O(kN^{3}+k^{2}N^{2})$ via fast Hankel SVD.

For the implementation comparison we consider the series $F_{N} =
(f_{1},\dots,f_{N})$ of the form
\begin{equation}
    \label{eqn:hseries}
    f_{n} =
    \begin{cases}
        \sin\left(\frac{2\pi}{10}n\right) + 0.01\varepsilon_{n}, \quad & n < Q, \\
        \sin\left(\frac{2\pi}{10.5}n\right) + 0.01\varepsilon_{n}, \quad & n \geq Q,
    \end{cases}
\end{equation}
where $\varepsilon_{n}$ denotes the uncorrelated white noise.

\begin{figure}[htbp]
    \centering
    \includegraphics[width=\columnwidth]{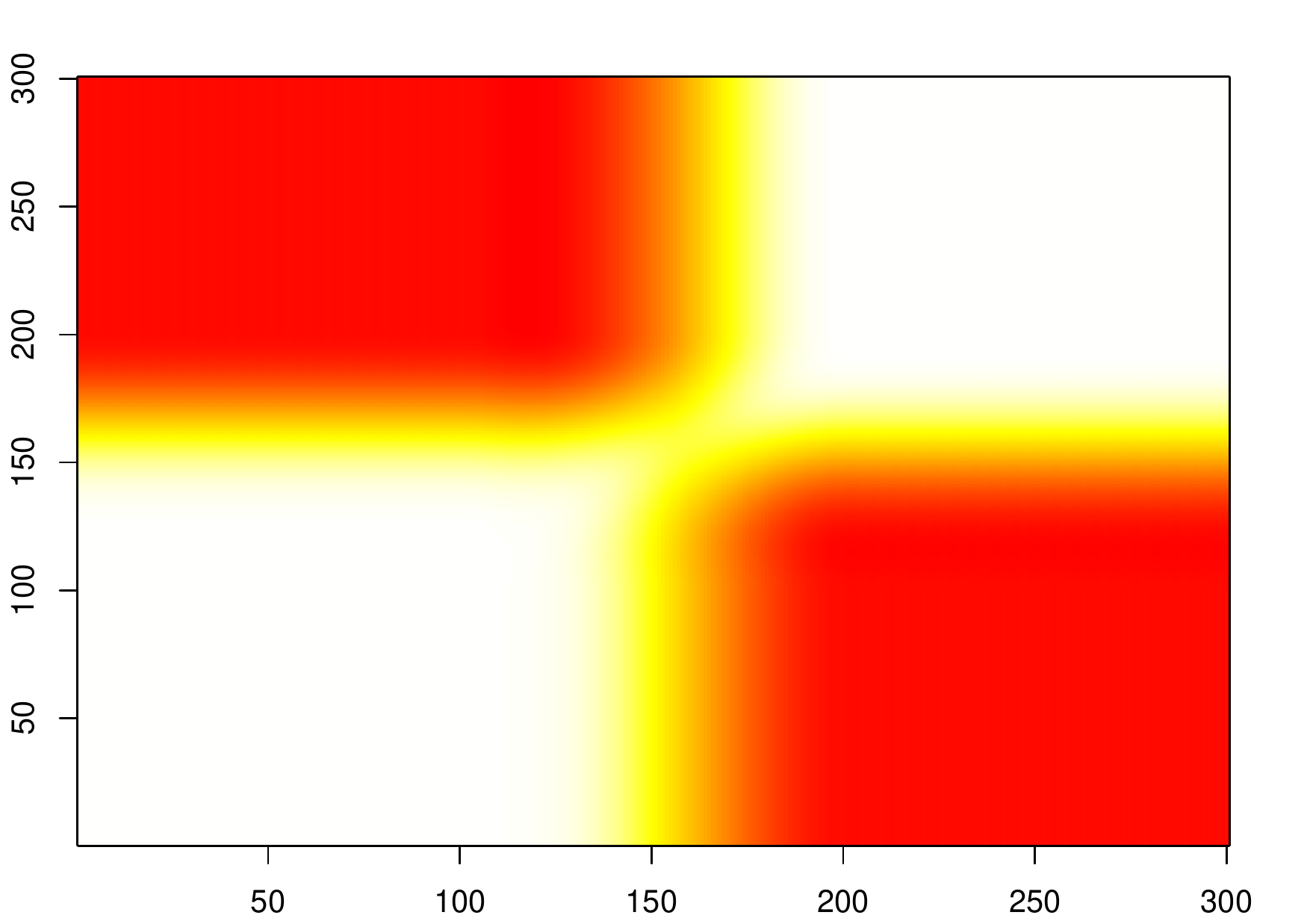}
    \caption{H-matrix for the series \eqref{eqn:hseries}.}
    \label{fig:hmatr}
\end{figure}

The typical H-matrix for such series is shown on the figure \ref{fig:hmatr}
(compare with figure 3.3 in \cite{Golyandina01}). The parameters used for this
matrix are $N=400$, $Q=200$, $B = 100$, $T = 100$, $L = 50$ and $I =
\left\{1,2\right\}$.

To save some computational time we compared not the worst case complexities, but
some intermediate situation. The parameters used were $Q=N/2$, $B = T = N/4$, $L
= B/2$ and $I = \left\{1,2\right\}$. For full SVD \verb@DGESDD@ routine from
\verb@ATLAS@ library \cite{Whaley04} was used (as it was shown below, it turned
out to be the fastest full SVD implementation available for our platform).

\begin{figure}[htbp]
    \centering
    \includegraphics[width=\columnwidth]{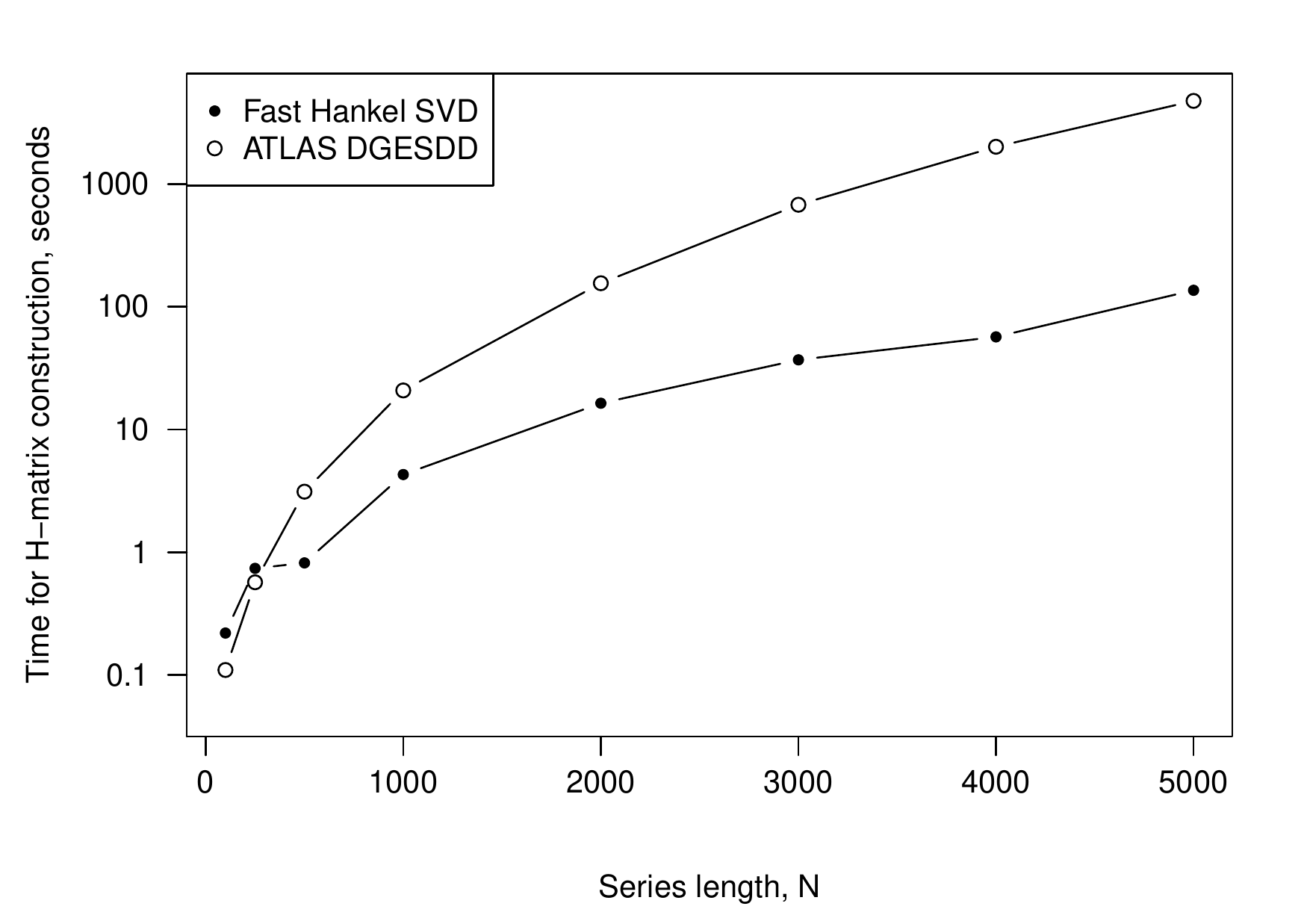}
    \caption{Structure Change Detection Comparison}
    \label{fig:strdet}
\end{figure}

The obtained results are shown on the figure \ref{fig:strdet}. As before,
logarithmic scale was used for the y-axis. Notice that the implementation of
series structure change detection using fast Hankel SVD is more than $10$ times
faster even on series of moderate length.

\section{Real-World Time Series}
\label{sec:example}
Fast SSA implementation allows us to use large window sizes during the
decomposition, which was not available before. This is a crucial point in the
achieving of asymptotic separability between different components of the series
\cite{Golyandina01}.

\subsection{Trend Extraction for HadCET dataset}
Hadley Centre Central England Temperature (HadCET) dataset \cite{Parker92} is
the longest instrumental record of temperature in the world. The mean daily data
begins in $1772$ and is representative of a roughly triangular area of the United
Kingdom enclosed by Lancashire, London and Bristol. Total series length is
$86867$ (up to October, 2009).

\begin{figure}[Htbp]
    \centering
    \includegraphics[width=\columnwidth]{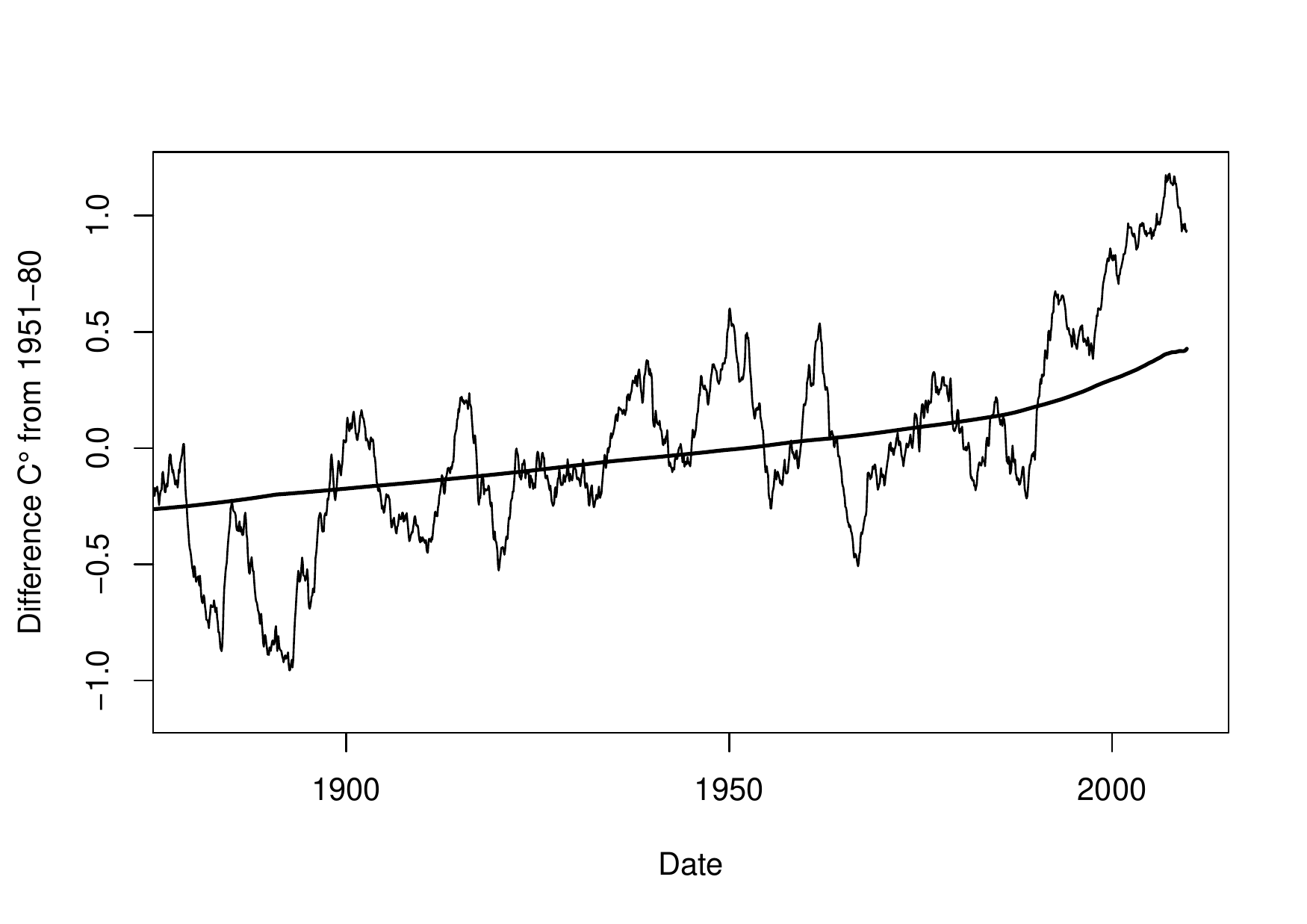}
    \caption{Hadley CET Series and Extracted Trend}
    \label{fig:cet}
\end{figure}

We applied SSA with the window length of $43433$ to obtain the trend and few
seasonal components. The decomposition took $22.5$ seconds and $50$ eigentriples
were computed. First eigentriple was obtained as a representative for the trend
of the series. The resulting trend on top of 5-year running mean temperature
anomaly relative to $1951$-$1980$ is shown on figure \ref{fig:cet} (compare with
figure~$1$ in \cite{Hansen06}).

\subsection{Quebec Birth Rate Series}
The series obtained from \cite{Hipel94} contains the number of daily births in
Quebec, Canada, from January 1, 1977 to December 31, 1990. They were discussed
in \cite{Golyandina01} and it was shown that in addition to a smooth trend, two
cycles of different ranges: the one-year periodicity and the one-week
periodicity can be extracted from the series.

However, the extraction of the year periodicity was incomplete since the series
turned out to be too smooth and do not uncover the complicated structure of the
year periodicity (for example, it does not contain the birth rate drop during
the Christmas). Total series length is $5113$, which suggests the window length
of $2556$ to be used to achieve the optimal separability. This seems to have
been computationaly expensive at the time of book written (even nowadays the
full SVD of such matrix using \verb@LAPACK@ takes several minutes and requires
decent amount of RAM).

\begin{figure}[htbp]
    \centering
    \includegraphics[width=\columnwidth]{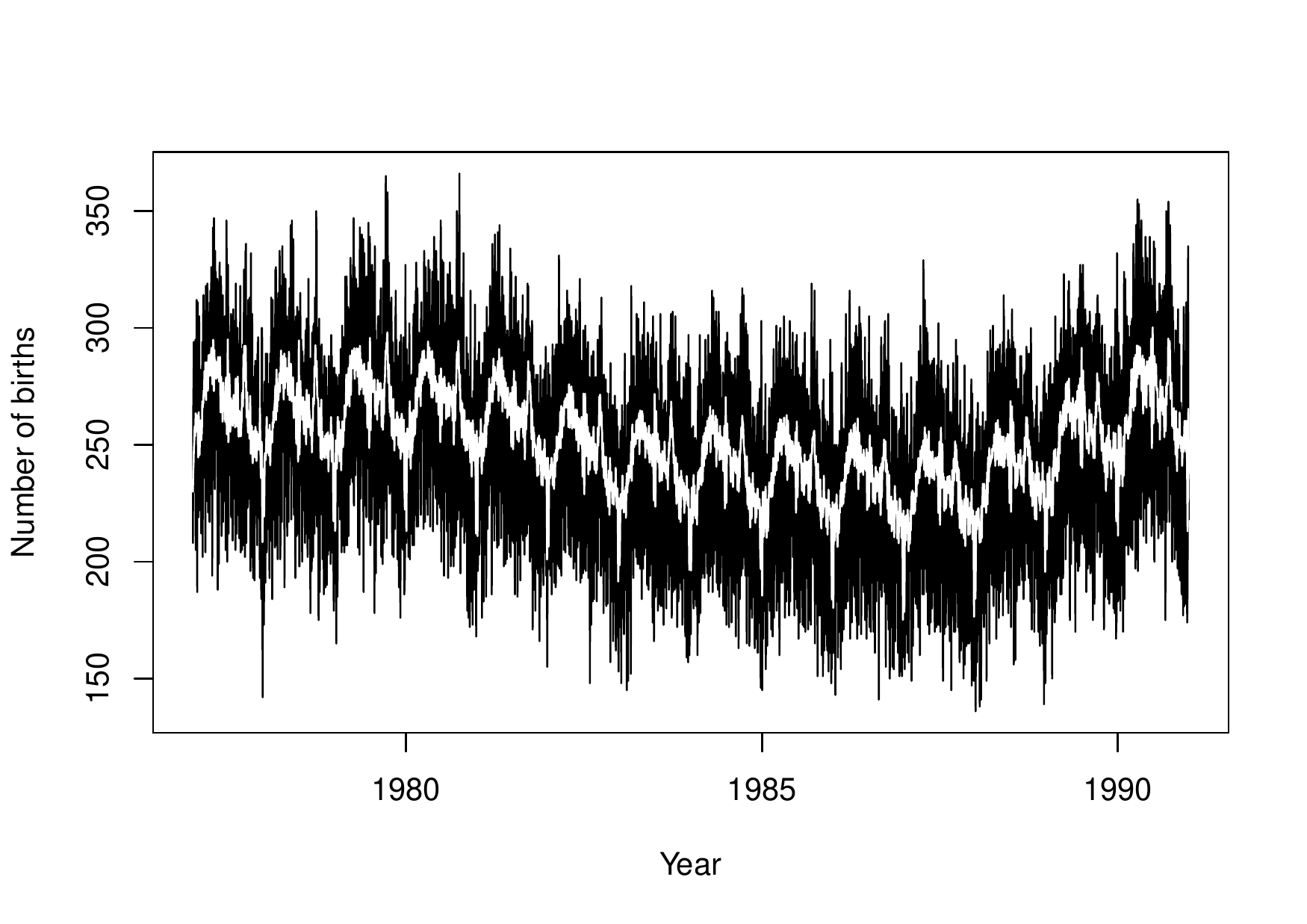}
    \caption{Quebec Birth Rate: Trend and Year Periodicity}
    \label{fig:births}
\end{figure}

We apply the sequential SSA approach: first we did the decomposition using
window size of $365$. This allows us to extract the trend and week periodicity
in the same way as described in \cite{Golyandina01}. After that we performed the
decomposition of the residual series using the window length of $2556$. $100$
eigentriples were computed in $3.2$ seconds. The decomposition cleanly revealed
the complicated structure of the year periodicity: we identified many components
corresponding to the periods of year, half-year, one third of the year, etc. In
the end, the eigentriples $1-58$ excluding $22-24$ were used during the
reconstruction. The end result containing trend and year periodicity is shown on
the figure~\ref{fig:births} (compare with figure~$1.9$ in
\cite{Golyandina01}). It explains the complicated series structure much better.

\section{Conclusion}
The stages of the basic SSA algorithm were reviewed and their computational and
space complexities were discussed. It was shown that the worst-case
computational complexity of $O(N^{3})$ can be reduced to $O(kN\log{N} +
k^{2}N)$, where $k$ is the number of desired eigentriples.

The implementation comparison was presented as well showing the superiority of
the outlined algorithms in the terms of running time. This speed improvement
enables the use of the SSA algorithm for the decomposition of quite long time
series using the optimal for asymptotic separability window size.

\section*{Acknowledgements}
Author would like to thank N.E.~Golyandina and anonymous reviewer for thoughtful
suggestions and comments that let to significantly improved presentation of the
article.

\bibliography{compssa}

\end{document}